\documentclass[11pt]{amsart}

\usepackage{amsmath}
\usepackage{amssymb}
\usepackage{amsfonts}
\usepackage{amsthm}
\usepackage{enumerate}

\usepackage{graphicx}
\usepackage{verbatim}

\def\Enu{\exp L^\nu}
\def\dyad{{\text{dyad}}}
\def\mart{{\text{dyad}}}
\def\intslash{\rlap{\kern  .32em $\mspace {.5mu}\backslash$ }\int}
\def\qsl{{\rlap{\kern  .32em $\mspace {.5mu}\backslash$ }\int_{Q_x}}}

\def\emph#1{{\it #1 }}

\def\ga{\gamma}

\def\cf{{\it cf}}

\def\lc{\lesssim}

\def\eps{\varepsilon}

\def\la{\lambda}

\def\bbD{{\mathbb {D}}}
\def\bbE{{\mathbb {E}}}

\def\bbN{{\mathbb {N}}}

\def\bbR{{\mathbb {R}}}

\def\bbT{{\mathbb {T}}}

\def\bbZ{{\mathbb {Z}}}

\def\cC{{\mathcal {C}}}

\def\cE{{\mathcal {E}}}

\def\cG{{\mathcal {G}}}

\def\cL{{\mathcal {L}}}

\def\cN{{\mathcal {N}}}

 at 10 true pt

\def\be#1{\begin{equation}\label{#1}}
\def\endeq{\end{equation}}
\def\endal{\end{align}}
\def\bas{\begin{align*}}
\def\eas{\end{align*}}
\def\bi{\begin{itemize}}
\def\ei{\end{itemize}}

\def\eps{\varepsilon}

\def\emph#1{{\it #1}}
\def\textbf#1{{\bf #1}}

\theoremstyle{plain}
   \newtheorem{theorem}[subsection]{Theorem}

   \newtheorem{proposition}[subsection]{Proposition}
   
   \newtheorem{lemma}[subsection]{Lemma}
   \newtheorem{corollary}[subsection]{Corollary}

\theoremstyle{remark}

\theoremstyle{definition}
   \newtheorem{definition}[subsection]{Definition}

\begin{document}

\title{Low regularity classes and entropy numbers}

\dedicatory
{In memory of Eduard Belinsky (1947 -- 2004)}

\author{ Andreas Seeger\  \  \   \ \ \ Walter Trebels}

\address{A. Seeger, Department of Mathematics, University of Wisconsin,
480 Lincoln Drive,
Madison, WI, 53706, USA}
\email{seeger@math.wisc.edu}

\address{W. Trebels, Fachbereich Mathematik, Technische Universit\"at
Darmstadt, Schlo\ss\-gar\-ten\-stra\ss e 7, 64289 Darmstadt,
Germany}
\email{trebels@mathematik.tu-darmstadt.de}

\subjclass{42B15}

\thanks{A.S.
 was supported in part by
an NSF grant.}

\date{March 1, 2007}

\begin{abstract}We note a sharp embedding of the Besov space 
$B^\infty_{0,q}(\bbT)$ into exponential classes and prove entropy
estimates  for the compact embedding of subclasses with logarithmic
smoothness,  considered by Kashin and Temlyakov.
\end{abstract}

\maketitle

\section{Introduction}

We consider spaces of functions with low regularity
and their embedding properties 
with respect to the exponential classes 
$\exp(L^\nu)$. For simplicity we work with functions 
 on the torus $\bbT=\bbR/\bbZ$ 
(identified with $1$-periodic functions on $\bbR$).
We use the following characterization of the Luxemburg norm in
$\Enu(\bbT)$, found
for example in \cite{Tr93}. For $\nu>0$ set
\begin{equation}
\|f\|_{\Enu(\bbT)}= \sup_{1\le p<\infty} p^{-1/\nu} \|f\|_{L^p(\bbT)};
\end{equation}
this norm will be used in what follows.

We consider the Besov spaces 
$B^\infty_{0,q}$, defined
via dyadic  decompositions as follows.
Let $\Phi\equiv\phi_0$ be an even $C^\infty$ function on $\bbR$
 with the property that $\Phi(s)=1$ for $|s|\le 1$ and $\Phi$ is supported
in $(-2,2)$. For $k\ge 1$ set $\phi_k(s)=\Phi(2^{-k} s)-\Phi(2^{-k+1}s)$ and,
for $k=0,1,2,\dots$
$$L_k f(x)\equiv \phi_k(D) f(x) = \sum_{n} \phi_k(n)\widehat {f_n} 
e^{2\pi i nx}.$$
Then 
$B^\infty_{0,q}$ is defined as the space of distributions for which
$$
\|f\|_{B^\infty_{0,q}}= \Big(\sum_{k=0}^\infty \big\|L_k
f\big\|_\infty^q \Big)^{1/q}
$$
is finite. It is well known that the class of functions defined in this way
does  not depend on the specific choice of $\Phi$.

The space $B^\infty_{0,q}$ consists of locally integrable functions 
if and only if $q\le 2$ (see \cite{har}, p. 112) 
 and it follows easily from the
definition that it embeds into $L^\infty$ if $q\le 1$.
We shall show for the interesting range $1<q\le 2$ a sharp embedding
 result 
involving the exponential classes.

\begin{theorem}\label{besov}
Let $1<q\le 2$. Then the space  
$B_{0,q}^{\infty}$
is continuously embedded in $\exp L^{q'}$, $q'=q/(q-1)$.
\end{theorem}

This can be read as a statement about the {\it growth envelope}
of the space $B_{0,q}^{\infty}$, defined by
\begin{equation}\label{grenv}
\cE_q(t)=\sup\{ f^*(t):  \|f\|_{B^\infty_{0,q}}\le 1\};
\end{equation}
here $f^*$ is the nonincreasing rearrangement of $f$.
It is shown in
Corollary 2.3 of \cite{ek}  that
$\|f\|_{\exp L^{q'}}\approx \sup_{t>0} f^*(t) \log^{-1/q'}(e/t)$ so  that
Theorem \ref{besov} immediately implies an upper bound $C|\log t|^{1/q'}$
for $\cE_q(t)$ when $t$ is small. The corresponding lower bound is proved
in \cite{har}, Prop. 8.24 (there also the  nonoptimal upper bound 
$C|\log t|$ is derived). Thus we get
\begin{corollary} For $1\le q\le 2$,
$$\cE_q(t)\approx |\log t|^{1/q'},\quad |t|\le
1/2.
$$
\end{corollary}

We shall now consider subclasses
$LG^\ga(\bbT)$ of
$B^\infty_{0,2}$ 
which are compactly embedded in Lebesgue and exponential classes; 
these were introduced by  Kashin and Temlyakov \cite{kt}.
For $\gamma>1/2$ the class $LG^\ga(\bbT)$  is defined as the class of 
$L^1(\bbT)$ functions
for which $\|L_k f\|_\infty = O((1+k)^{-\gamma})$ and we set
$$\|f\|_{LG^\ga(\bbT)} =\sup_{k\ge 0} (1+k)^{\gamma} \|L_k f\|_\infty.$$
Clearly, for $\gamma>1$ the class 
$LG^\ga(\bbT)$ is embedded in $L^\infty$ and if $1/2<\gamma\le 1$  
then $LG^\ga(\bbT)$ is embedded in $\Enu(\bbT)$ for $\nu<(1-\gamma)^{-1}$,
by Theorem \ref{besov}. We are interested in the compactness properties
of  this embedding and some related quantitative statements.

We recall that given a Banach space $X$ and a subspace $Y\subset X$ 
one defines the
{\it $n$th
entropy number} $e_n(Y;X)$
as the infimum over all numbers $\eps>0$ for which there are $2^{n-1}$
balls  of radius $\eps$ in $X$ which cover the  unit ball $\{y\in Y:
\|y\|_Y \le 1\}$ embedded  in $X$.
It is easy to see that  the embedding of $Y$ in $X$ is a compact operator 
if and only if $\lim_{n\to\infty} e_n(Y;X)=0$.

For $\gamma>1$ the embedding of $LG^\ga(\bbT)$ 
into $L^\infty$ is compact and 
Kashin and Temlyakov \cite{kt} 
 determined sharp bounds for the entropy numbers for the 
embedding into $L^\infty$ and $L^p$, $p<\infty$;  they showed
that for $n\ge 2$ and $\gamma>1$
\begin{equation}\label{KTresult}
e_n(LG^\gamma, L^p)\approx\begin{cases}
(\log n)^{1/2-\gamma},  &1\le p<\infty,
\\
(\log n)^{1-\gamma}, &p=\infty.
\end{cases}
\end{equation}

We note that the restriction $\gamma>1$ in \cite{kt} is only used to 
ensure 
the imbedding into $L^\infty$; indeed it is implicitly in \cite{kt} that
for $p<\infty$ the $L^p$ result \eqref{KTresult} holds for all
 $\gamma>1/2$. 
The hard part in the Kashin-Temlyakov result are the lower bounds.
 The $L^p$ lower bound  is derived using Littlewood-Paley theory from
lower  bounds for classes of trigonometric polynomials 
in  \cite{kash}. The $L^\infty$ bounds require fine estimates 
for certain Riesz products (\cf. Theorem 2.3 in \cite{kt}).

It is desirable to explain  the jump in the exponent 
that occurs in \eqref{KTresult} when $p\to \infty$. To achieve this
Belinsky and Trebels (\cite{bt}, Theorem 5.3) studied 
the entropy numbers 
$e_n(LG^\gamma,\Enu)$
for the natural embedding into the exponential classes; 
they obtained the equivalence
$e_n\approx (\log n)^{1/2-\gamma}$ for $\nu\le 1$.
For $\nu\ge 2$ they obtained an almost sharp result, namely that
$e_n$ is essentially $(\log n)^{1-\gamma-1/\nu}$, 
albeit with a  loss of $(\log\log n)^{1/\nu}$ for the upper bound. 
A more substantial gap between lower and upper bounds 
 remained  for $1\le \nu<2$. In \cite{bt} it was
also noticed that this gap could be closed if Pichorides conjecture
\cite{pich}  on the constant in the reverse 
 Littlewood-Paley inequality  were proved;  this however is still  an 
open problem. Nevertheless we shall use this insight to close the 
gap in \cite{bt}.


\begin{theorem}\label{mainresult} The embedding
$LG^\gamma(\bbT)\to \Enu(\bbT)$ is compact if either
$\gamma>1/2$, $\nu<2$,  or $\nu\ge 2$, $\gamma>1-\nu^{-1}$, and there
are the following upper and lower bounds for the entropy numbers.

(i) For $\gamma>1/2$, and $\nu<2$,
\begin{equation} \label{equiv1}
e_n(LG^\gamma, \Enu) \approx
(\log n)^{1/2-\gamma}.
\end{equation}

(ii) For $\nu\ge 2$ and $\gamma>1-\nu^{-1}$,
\begin{equation}\label{equiv2}
e_n(LG^\gamma, \Enu) \approx
(\log n)^{1-\gamma-1/\nu}.
\end{equation}
\end{theorem}

The lower bounds are known; for $\nu\le 2$ they follow 
immediately from \eqref{KTresult}. It was pointed out in \cite{bt} that
for $\nu>2$ the lower bounds   follow from the $L^\infty$ lower bound in 
\eqref{KTresult} and $L^\infty\to \exp(L^\nu)$ Nikolskii inequalities 
for trigonometric polynomials.

We thus are left to  establish the upper bounds for the entropy numbers.
The idea here is to embed the classes $LG^\gamma$ into slightly larger
classes $LG_\dyad^\gamma$  which contain discontinuous functions but 
satisfy the same entropy estimates with respect to the exponential classes.
Instead of the Pichorides conjecture we shall then use the well known
bounds  for a martingale analogue, due to Chang, Wilson and Wolff
\cite{cww}. 
This philosophy also applies to the proof of Theorem \ref{besov}; it has
been  used in other papers, among them
\cite{jkrw}, \cite{jsw}, \cite{ghs} (see also references contained in
these  papers).

{\it Notation.}
If $\, X, \, Y $ are normed linear spaces we use the notation $\, Y
\hookrightarrow X$ to indicate that $ \, Y \subset X$ and the embedding
is continuous.

{\it This paper.} The proof of  Theorem \ref{besov} is given 
 in \S\ref{embedd}, and the proof of  Theorem \ref{mainresult} in 
\S \ref{KashTe}.

\section{Embedding into the exponential classes}\label{embedd}

We shall work with dyadic versions of the Besov spaces where the 
Little\-wood-Paley operators $L_k$ are replaced 
by martingale difference operators.
Let $k$ be a nonnegative integer.
For a function on $[0,1]$ we define the conditional expectation operator
$$
\mathbb  E_k f(x) = 2^k \int_{(m-1) 2^{-k}}^{m2^{-k}} f(t) dt ,\quad
(m-1) 2^{-k}\le x<m2^{-k} ,\quad  m=1,\dots, 2^k,
$$
and define
\begin{align*}
\mathbb  D_k f(x)&=\mathbb  E_{k} f(x)-\mathbb  E_{k-1} f(x), \quad k\ge 1,\\
\mathbb  D_0 f(x)&=\mathbb  E_{0} f(x);
\end{align*}
clearly both $\bbE_k f$ and $\bbD_k f$ define $1$-periodic functions and
can  be viewed as functions on $\bbT$.
Note that the functions $\bbD_k f$ are piecewise constant and (typically)
discontinuous at $m2^{-k}$, $m=0,\dots, 2^k-1$. We also observe that 
$ f = \sum_{k\ge 0}^{} \mathbb  D_k f \, $ almost everywhere for $\, f
\in L^1 .$

\begin{definition}
Let $1\le q\le 2$.
The dyadic Besov-type spaces $\ell^{q}(B^\infty_\mart)$ consists of all 
$f\in L^1(\bbT)$ for which the sequence 
$\{\|\bbD_k f\|_\infty\}_{k=0}^\infty$ belongs to $\ell^{q}$; the norm is
given by
$$
\|f\|_{\ell^{q}(B^\infty_\mart)}=
\Big(\sum_{k=0}^\infty\|\bbD_k f\|_\infty^q\Big)^{1/q}.
$$
\end{definition}

\begin{proposition} \label{dy}
Let $1\le q\le 2$. Then
$$
B^\infty_{0,q}\hookrightarrow \ell^{q}(B^\infty_\mart)\, .
$$
\end{proposition}

This is easily reduced to the following  estimate on compositions of
the difference operators with the convolutions $\phi(D/\la)$ for large
$\la$. 

\begin{lemma} \label{contversdyad}
Let $\la\ge 1$ and $k\ge 0$. Let
$\psi\in C^\infty$ be even,  with support in $(-2,-1/2)\cup(1/2,2)$ and
let $\cL_\la=\psi(\la^{-1}D)$. Then
\begin{align}
\label{Ekphila}
\big\|\bbE _{k} \cL_\la\big\|_{L^\infty\to L^\infty} &\le C
\min\{ \la^{-1} 2^k,1 \}, \quad k\ge 0,
\\
\label{Dkphila}
\big\|\bbD_k \cL_\la\big\|_{L^\infty\to L^\infty} &\le C
\min\{ \la^{-1} 2^k, \la 2^{-k}\}, \quad k\ge 1.
\end{align}
\end{lemma}

\begin{proof}
Use the notation $ \psi _{-1}(s)=(2\pi i s)^{-1}\psi(s), \psi _1 (s) =s \,
\psi(s)$  and observe  that $\psi, \psi _{-1}, \psi _1$ are
$C^\infty$--functions with compact support away from the origin so that 
 by standard $\widehat{L^1}$-theory
 the sequences  $\ell\mapsto \psi(\la^{-1}\ell)$,   
$\psi _{-1}(\la^{-1}\ell)$,  $\psi _1(\la^{-1}\ell) $ 
define the Fourier coefficients 
of $L^1(\bbT)$ functions, with $L^1$ norms uniformly in $\la \, .$ Therefore, 
\begin{equation}\label{Ekphila*}
\|\psi(\la^{-1}D) f\|_\infty
+
\| \psi _{-1}(\la^{-1}D )f\|_\infty + \|\psi _1(\la^{-1}D) f\|_\infty 
\le C\|f\|_\infty \, .
\end{equation}
In particular it is clear that
$\|\bbE_k \cL_\la\|_{L^\infty\to L^\infty}=O(1)$.

Now fix $k$ so that $2^k<\la$ and let  $x_{m,k}=m 2^{-k}$. Then
for $x\in [x_{m,k}, x_{m+1,k})$,
\begin{align*}
\bbE_{k} \cL_\la f(x) 
&= 2^k\int_{x_{m,k}}^{x_{m+1,k}} \left(  
\sum_{\ell\in\bbZ}\psi(\la^{-1}\ell)\int_0^1
e^{-2\pi i \ell y} f(y)\,  dy \, e^{2\pi i \ell x} \right)  dx
\\
&=2^k\sum_{\ell\in\bbZ}\psi(\la^{-1}\ell)\int_0^1
\frac{e^{2\pi i \ell(x_{m+1,k}-y)}- e^{2\pi i \ell(x_{m,k}-y)}}{2\pi i
\, \ell}  f(y) \, dy
\\
&=  2^k\la^{-1}\Big( \psi _{-1}(D/\la)f(x_{m+1,k})
- \psi _{-1} (D/\la)f(x_{m,k})\Big)
\end{align*}
and \eqref{Ekphila} follows by \eqref{Ekphila*}.

Inequality \eqref{Dkphila}
for $2^k<\la$ is an immediate consequence and it remains  to consider
the case $2^k\ge \la$.
Fix $x$, then $\bbE_k \cL_\la f(x)$ is the average of $\cL_\la f$ 
over an interval of length $2^{-k}$ containing $x$. Thus, by the 
mean value theorem applied to $\bbE_k \cL_\la  f(x)$ and 
$\bbE_{k-1} \cL_\la  f(x)$, we can write for $k\ge 1$
$$
\bbD_k \cL_\la f(x) = \cL_\la f(x')-\cL_\la f(x'')=
(\cL_\la f)'(\tilde x) (x'-x'')
$$
where $x', x'', \tilde x$ have distance at most $2^{-k+1}$ from $x$.
Now $(\cL_\la f)' = \la \psi _1(D/\la) f$ 
and thus
$$ 
\|\bbD_k \cL_\la f\|_\infty \le 2^{1-k} \|(\cL_\la f)'\|_\infty  
\le C\la 2^{-k}  \|f \|_\infty.
$$
\end{proof}

\begin{proof}[Proof of Proposition \ref{dy}]
Let $\Psi_0$ be a $C^\infty$ function supported in $(-4,4)$
which satisfies $\Psi_0(s)=1$ in $(-2,2)$ and let
$\Psi_n= \Psi(2^{-n}\cdot)$ where $\Psi $ is supported in
$(-8, -1/8)\cup(1/8,8)$ so that $\Psi(s)=1$ for $|s|\in (1/2,4)$.
Then $\Psi_n\phi_n=\phi_n$ for all $n$, so that $\Psi_n(D) L_n=L_n$, 
and we can
write
\begin{align*}
\big\|\bbD_k f\big\|_\infty&= \Big\|
\bbD_k \sum_{n=0}^\infty \Psi_n(D) L_n f\Big\|_\infty
\le \sum_{n=0}^\infty \big\|\bbD_k \Psi_n(D)
\big\|_{L^\infty\to L^\infty}
\big\|
L_n f\big\|_\infty
\\&
\le C\sum_{n=0}^\infty 
2^{-|k-n|}\|L_n f\|_\infty
\end{align*}
and therefore 
\begin{equation*}
\|f\|_{\ell^{q}(B^\infty_\mart)}
\le C \sum_{m=0}^\infty 2^{-m}\Big\|\big\{ \|L_{k+m} f\|_\infty 
\big\}_{k=-m}^\infty\Big\|_{\ell^{q}} \le C'
\|f\|_{B^\infty_{0,q}}.
\end{equation*}
\end{proof}

We now introduce the  square-function and the maximal function
$$
\frak S (f):=\Big(\sum_{k\ge 0}|\mathbb  D_k f(x)|^2\Big)^{1/2} \,
,\qquad \frak M _0 (f):= \sup_{k\ge 0} |\mathbb  E_k f(x)-\mathbb
E_0 f(x)|\, ,
$$
resp., and recall the following deep
``good $\lambda$ inequality'' due to Chang, Wilson and Wolff (Corollary
3.1
in \cite{cww}): There are absolute constants $c$ and $C$
so that for all  $\lambda>0$, $0<\varepsilon <1$,
\begin{multline}\label{glaineq}
\text{meas}\big(\big\{x: \frak M _0 (f)(x)
>2\lambda,\, \frak S(f)< \varepsilon\lambda\big\}\big)
\\ \leq C\, {\rm exp}(-\frac{c}{\varepsilon^2})
\text{meas}\big(\big\{x:\sup_{k\ge 0} |\mathbb  E_k f(x)|>
\lambda\big\}\big).
\end{multline}
It is standard that this implies the inequality
\begin{equation}\label{revLP}
\|f\|_p \le C \sqrt p \,\|\frak S(f)\|_p
\end{equation}   
for all $p\ge 2$, and some absolute constant $C\ge 1$.
Indeed, if we integrate out the $L^p$ norms using the distribution
function, where we observe that
$$
\{ x : \frak M _0 (f)(x) > 2\lambda \} \subset \{ x : \frak M _0 (f)
> 2\lambda ,\, \frak S (f)< \varepsilon\lambda\} \cup
\{ x : \frak S (f) \ge \varepsilon \lambda \}\, ,
$$
we obtain
\begin{equation*}
\big\|\sup_k |\mathbb E_k f|\big\|_p\le
\|\mathbb E_0 f\|_p + 2 C^{1/p} e^{-c\varepsilon^{-2}p^{-1}}
\|\sup_k|\mathbb E_k f|\big\|_p + 2\varepsilon^{-1}\big\|\frak S
(f) \big\|_p.
\end{equation*}
Now we  choose $\varepsilon= a p^{-1/2}$ with $a$ so small  that
$2 C e^{-ca^{-2}}=1/2.$  Since
 $\mathbb D_0=\mathbb E_0$ is incorporated in the definition of the
square-function, $\, |f(x)| \le \sup_k|\mathbb E_k g|(x)  $ a.e.,
the  asserted bound \eqref{revLP} follows.

The following interpolation result is a quick consequence of  
\eqref{revLP}.


\begin{lemma} \label{interpol} There is a constant $C$ so that for 
$1\le s\le 2$, $s'=s/(s-1)$, $2\le p<\infty$, and all sequences $\{f_k\}$ of $L^p(\bbT)$ functions,
\begin{equation*}
\Big\|\sum_{k=0}^\infty \bbD_k f_k \Big\|_{L^p(\bbT)} \le C p^{1/s'}
\Big(\sum_{k=0}^\infty \|f_k\|_{L^p(\bbT)}^s\Big)^{1/s}.
\end{equation*}
\end{lemma}

\begin{proof}
The statement is trivial for $s=1$, because of the uniform $L^p$ bounds
for  the operators $\bbD_k$.
We thus only need to prove the statement for $s=2$ since
 then the general case follows by complex interpolation.
By a straightforward limiting argument we
may assume that $f_k=0$ for all but finitely many $k$.

We use that $\bbD_k\bbD_l=0$ if $k\neq l$, and 
define $g=\sum \bbD_k f_k$.
Then by \eqref{revLP}
\begin{align*}
\|g\|_p= \Big\|\sum_l \bbD_l g\Big\|_p
\le C \sqrt p \Big\|\Big(\sum_l|\bbD_l g|^2\Big)^{1/2}\Big\|_p,
\end{align*}
and since $p\ge 2$ we can use  Minkowski's inequality to bound this by
\begin{align*}
C \sqrt p       \Big(\sum_l\|\bbD_l g\|_p^2\Big)^{1/2}
=C \sqrt p       \Big(\sum_l\|\bbD_l f_l\|_p^2\Big)^{1/2}
\le C' \sqrt p       \Big(\sum_l\|f_l\|_p^2\Big)^{1/2}.
\end{align*}\end{proof}

Theorem \ref{besov} is an immediate
consequence  of Proposition \ref{dy} and the following imbedding result
which is  based on \eqref{revLP} (or rather the case $s=2$ of 
Lemma \ref{interpol}).
\begin{proposition}\label{expimb}
Let $1\le q\le 2.$ Then 
\begin{equation*} 
 \ell^q(B^\infty_\mart)\hookrightarrow \exp L^{q'}.
\end{equation*}
\end{proposition}

\begin{proof} We modify an argument from \cite{bt}
(which was based there on the Pichorides  conjecture).
Fix $f\in \ell^q(B^\infty_\mart)$ and let $n\to k(n,f)$ be a bijection
of $\bbN\cup\{0\}$ so that the sequence $n\to \|\bbD_{k(n,f)} f\|_\infty$
is nonincreasing (in other words, we form the nonincreasing rearrangement
of the sequence $\{\| \bbD_k f \| \}$).

For $p\ge 2$ we need to estimate $p^{-1/q'}\|f\|_p$. Thus fix $p>2$ 
and let  $N\in \bbN$ so that $p\le N<p+1$.
We then split
$$
f=\sum_{n=0}^N \bbD_{k(n,f)} f + 
\sum_{n=N+1}^\infty
\bbD_{k(n,f)} f  := I_Nf + II_N f.
$$
By H\"older's inequality
\begin{align}
\|I_Nf\|_p &\le \sum_{n=0}^N \| \bbD_{k(n,f)} f\|_p
\le \sum_{n=0}^N \| \bbD_{k(n,f)} f\|_\infty
\notag
\\
\label{Ibd}
&\le (N+1)^{1/q'} 
\Big(\sum_{n=0}^N \| \bbD_{k(n,f)} f\|_\infty^q\Big)^{1/q} \le 
Cp^{1/q'} \|f\|_{\ell^q(B^\infty_\mart)}.
\end{align}

For the second term we get a bound in terms of the 
Lorentz-Besov  type space 
$\ell^{q,2}(B^\infty_\mart)$ defined similarly as
$\ell^{q}(B^\infty_\mart)$, but with the sequence space $\ell^q$ replaced 
by the Lorentz variant $\ell^{q,2}$. Since $\ell^q\subset \ell^{q,2}$ for
$q\le 2$; this is a better estimate.  Note that 
\begin{equation}\label{lorentzdef} 
\Big\| \{ \bbD_k f\}_{k=0}^\infty \Big\|_{\ell^{q,2}} \approx
\Big(\sum_{n=0}^\infty \big[ n^{1/q} \|\bbD_{k(n,f)} f\|_\infty\big]^2 
n^{-1} \Big)^{1/2}.
\end{equation}

We now use the case $s=2$ of Lemma \ref{interpol} to obtain
\begin{align*}
\|II_N f\|_p &\le C p^{1/2} \Big(\sum_{n=N+1}^\infty \big\|\bbD_{k(n,f)} f
\big\|_p^2\Big)^{1/2}
\le C p^{1/2} \Big(\sum_{n=N+1}^\infty \big\|\bbD_{k(n,f)} f
\big\|_\infty^2\Big)^{1/2}
\\
&\le C p^{1/2} N^{-1/2+1/q'}
\Big(\sum_{n=N+1}^\infty n^{1-2/q'}\big\|\bbD_{k(n,f)} f
\big\|_\infty^2\Big)^{1/2},
\end{align*}
and, since $1-2/q'=2/q-1$ and $p\approx N$, we get from
\eqref{lorentzdef} 
\begin{equation} \label{IIbd}
p^{-1/q'}\|II_N f\|_p 
\le C \|f\|_{\ell^{q,2}(B^\infty_\mart)}
\le C' \|f\|_{\ell^{q}(B^\infty_\mart)}.
\end{equation}
Estimates \eqref{Ibd} and \eqref{IIbd} yield 
$$\|f\|_{\exp L^{q'}}\lc \|f\|_{\ell^{q}(B^\infty_\mart)}$$
and thus the assertion. 
\end{proof}


\section{Entropy numbers for the Kashin-Temlyakov classes}\label{KashTe}
We now give a proof of Theorem \ref{mainresult}.
As discussed in the introduction only the upper bounds have to be proved.
 It will be advantageous to define larger ``dyadic'' analogues of 
the $LG$ classes.
\begin{definition}
Let $\gamma>1/2$ and let $LG^\gamma_\dyad(\bbT)$ denote 
the class of $L^1(\bbT)$ 
 functions for which
$\|\mathbb  D_k f\|_\infty =O(k^{-\gamma})$ as $k\to \infty$. We set
$$
\|f\|_{LG^\gamma_\dyad}= \sup_{k\ge 0} \, (k+1)^\gamma
 \|\bbD_k f\|_\infty.
$$
\end{definition}

We note that the classes $LG^\gamma(\bbT)$ consist of continuous functions
provided that $\gamma>1$. This is not the case for 
the dyadic analogue $LG^\gamma_\dyad(\bbT)$ as even the building blocks 
$\bbD_k f$ are piecewise constant and typically discontinuous 
at $m2^{-k}$, $m=0,\dots, 2^k-1$. We prove the following embedding result.

\begin{lemma} \label{star}
For $\gamma>1/2$
$$ 
LG^\gamma(\bbT)\hookrightarrow LG^\gamma_\dyad(\bbT)\, .
$$
\end{lemma}

\begin{proof}
This follows easily from Lemma
\ref{contversdyad}. Indeed let $f\in LG^\gamma(\bbT)$,
so that $\|L_n f\|_\infty \lesssim \|f\|_{LG^\gamma} 
(1+n)^{-\gamma}$. As in \S\ref{embedd} we can write
$L_n= \Psi_n(D) L_n$ where the operator $\bbD_k\Psi_n(D)$ has
$L^\infty\to L^\infty$  operator  norm $O( 2^{-|k-n|})$.
Thus
\begin{align*}
\big\|\bbD_k f\big\|_\infty&= \Big\|
\bbD_k \sum_{n=0}^\infty \Psi_n(D) L_n f\Big\|_\infty
\le C\sum_{n=0}^\infty 2^{-|k-n|}
\big\|
L_n f\big\|_\infty
\\
&\le C_0 \sum_{n=0}^\infty  2^{-|k-n|} (1+n)^{-\ga} \|f\|_{LG^\ga}
\le C' (1+k)^{-\ga} \|f\|_{LG^\ga}.
\end{align*}
This  proves the assertion.
\end{proof}

We now state a crucial  approximation result which will be derived as a 
quick consequence of Lemma \ref{interpol}.

\begin{lemma} \label{appr}
Let $1/2<\gamma<1$ and $0<\nu<(1-\gamma)^{-1}$
or $\gamma\ge 1$ and $0<\nu<\infty$.  There is a constant 
$C=C(\gamma,\nu)$ 
so that for 
$M=1,2,\dots$
\begin{equation*}
\sup_{\|f\|_{LG^\gamma_\dyad} \le 1} \, \|f-\bbE_M f\|_{\exp L^\nu}
 \le C
\begin{cases} M^{1/2-\gamma},
\quad &\nu\le 2, \,\gamma>1/2,
\\
M^{1-1/\nu-\gamma}, \quad &\nu\ge 2, \,
\gamma> 1-\nu^{-1}.
\end{cases}
\end{equation*}
\end{lemma}
\begin{proof}
Consider $\, f \in LG^\gamma_\dyad ,\, \|f\|_{LG^\gamma_\dyad} \le 1,$ and 
 write 
$$
f-\bbE_M f= \sum_{k=M+1}^\infty \bbD_k \bbD_k f \, .
$$ 
By Lemma \ref{interpol} we have for $2\le p<\infty$, and $s\gamma>1$
\begin{align}
&p^{-1/\nu} \|f-\bbE_M f\|_p \le C p^{1/s'-1/\nu}\Big(\sum_{k=M+1}^\infty
\big\|
\bbD_k f\big\|_p^s\Big)^{1/s}
\notag
\\&\,\le C
p^{1/s'-1/\nu}
\Big(\sum_{k=M+1}^\infty
\big\|
\bbD_k f\big\|_\infty^s\Big)^{1/s}
\le C p^{1/s'-1/\nu}
\Big(\sum_{k=M+1}^\infty
(1+k)^{-s\gamma}\Big)^{1/s}
\notag
\\&\le \cC(s,\ga) p^{1-1/\nu-1/s}
 M^{1/s-\gamma}.
\notag
\end{align}
If $\nu\le 2$ then we may apply this bound for $s=2$, $\gamma>1/2$ and
get  the bound $\|f-\bbE_M f\|_{\exp L^\nu}=O(M^{1/2-\gamma})$. 
If $\nu >2$ we may apply it with $s=\nu/(\nu-1)\in (1,2)$, indeed 
 we have $s\gamma>1$ in view of our assumption $\gamma>1-\nu^{-1}$; the
result  is the asserted bound $\|f-\bbE_M f\|_{\exp
L^\nu}=O(M^{1-1/\nu-\ga})$. 
\end{proof}

We  apply a result of Lorentz \cite{l}, \cf.  Theorem 3.1 in  
\cite{lgm}, p. 492. Here one considers  a Banach space $X$ of functions,
a sequence $\cG=\{g_1,g_2,\dots\}$ of linearly independent functions whose
linear span is dense in $X$.  Set $X_0=0$, and let, for $n\ge 1$, $X_n$ 
be the linear span of $g_1,\dots, g_n$. Let 
$$
D_n(x)= \inf\{\|x-y\|: y\in X_n\}
$$ 
and let 
$\frak d=(\delta_0,\delta_1,\dots)$ be a nonincreasing sequence of
positive  numbers with $\lim_{n\to \infty}\delta_n=0$.
Let 
$$
A(\frak d)=\{x\in X : D_n(x)\le \delta_n, \,n=0,1,2,\dots\}
$$
be the approximation set  associated with $\frak d$, $\cG$.

Next let  $\cN_\eps(A(\frak d))$ denotes the minimal number 
of balls of radius $\eps$ needed to cover  $A(\frak d)$.
The following inequality for the natural logarithm of $\cN_\eps(A(\frak
d))$  is a special case of Lorentz' result.
\begin{equation} \label{lorentzestimate}
\log N_\eps(A(\frak d))\le 2n \log \big(\tfrac{18\delta_0}\eps\big), 
\quad \text{ if } \eps\ge \delta_n \,  .
\end{equation}

We apply 
\eqref{lorentzestimate}
 to  prove the dyadic analogue of the upper bound in
Theorem \ref{mainresult}.
\begin{proposition}\label{entrdyad}
The embedding
$LG^\gamma_\dyad(\bbT)\to \Enu(\bbT)$ is compact if
$\gamma>1/2$, $\nu<2$ or 
$\nu\ge 2$, $\gamma>1-\nu^{-1}$ and we have 
\begin{align}\label{entrdyadest1}
e_n(LG^\gamma_\dyad, \Enu) &\le C
(\log n)^{1/2-\gamma}, \quad \gamma>1/2,\quad \nu\le 2,
\\
\label{entrdyadest2}
e_n(LG^\gamma_\dyad, \Enu) &\le C
(\log n)^{1-\gamma-1/\nu}, \quad \gamma>1-1/\nu, \quad\nu\ge 2.
\end{align}
\end{proposition}

\begin{proof}
We set $X=\Enu$, and, for $n=2^M+j$, $j=0,1,\dots, 2^M-1$, let $g_n$ be
the  characteristic function of the interval $[j 2^{-M}, (j+1) 2^{-M})$. 
If $X_n$, $D_n(x)$ are defined as above then
we note that Lemma \ref{appr} says that for $f $ in the unit ball of 
$LG^\gamma_\dyad$ we have
$$
D_n(f)\le C_0 (\log (n+2))^{-a}
$$
where $a=\gamma-1/2$  if
$\gamma>1/2$ and $\nu\le 2$, and 
$a=\gamma+\nu^{-1}-1$ if 
$\nu\ge 2$ and $\gamma>1-1/\nu$.
We now note that  \eqref{lorentzestimate} implies
that
$$
e_{\widetilde n}(LG^\gamma_\dyad, \Enu) \le (\log (n+1))^{-a}
$$
if $\widetilde n> C n \log\log n$.
As $\log \widetilde n\approx \log n$ the asserted inequalities follow.
\end{proof}

\begin{proof}[Conclusion of the proof]
By Lemma \ref{star} 
we have
\begin{equation}\label{comp} 
e_{n}(LG^\gamma, \Enu) \le C 
e_{n}(LG^\gamma_\dyad, \Enu) 
\end{equation} 
and the assertion of the Theorem \ref{mainresult} follows 
from Proposition 
\ref{entrdyad}. 
\end{proof}

{\it Remark:} We note that  in the dyadic case, 
 there are also similar  lower bounds  matching
\eqref{entrdyadest1}, \eqref{entrdyadest2} 
for the
entropy numbers $e_n(LG^\gamma_\dyad, \Enu)$. These follow from
\eqref{comp} and 
 the known lower bounds for the entropy numbers for  $LG^\gamma$.

\end{document}